\documentclass[12pt,a4paper,reqno]{amsart}
\usepackage{amsmath,amsthm,amscd}
\usepackage{amsfonts,amssymb}
\allowdisplaybreaks

\newcommand{\be}{\begin{equation}}
\newcommand{\ef}{\end{equation}}
\chardef\bslash=`\\ 

\newtheorem*{thm*}{Theorem}

\theoremstyle{definition}

\newtheorem*{remark*}{Remarks}
\newtheorem*{defn*}{Definition}

\theoremstyle{remark}

\newcommand{\wt}{\widetilde}
\newcommand{\wh}{\widehat}

 \renewcommand{\sectionmark}[1]{}
\renewcommand{\Im}{\operatorname{Im}}

\newcommand{\ve}{\varepsilon}

\newcommand{\iy} {\infty}
\newcommand{\fc} {\frac}

\usepackage{amsfonts}
\newcommand{\field}[1]{\mathbb{#1}}

\newcommand{\dl}{\delta}
\newcommand{\D}{\field D}

\newcommand{\z}{\zeta}
\newcommand{\ov}{\overline}

\newcommand{\hC}{\wh{\field{C}}}
\newcommand{\C}{\field{C}}

\newcommand{\B}{\mathbf{B}}
\newcommand{\T}{\mathbf{T}}

\newcommand{\Belt} {\operatorname{Belt}}

\newcommand{\const}{\operatorname{const}}

\newcommand{\Om} {\Omega}
\newcommand{\vk} {\varkappa}
\newcommand{\x} {\mathbf x}
\renewcommand{\a} {\alpha}

\newcommand{\ld}{\lambda}

\begin{document}

\title{Quantitative theory of reflections across quasiconformal polygonal lines}
\author{Samuel L. Krushkal}

\begin{abstract} The paper continues the author's research in the problem of quantitative investigation of basic curvelinear quasiinvariants of quasiconformal curves. It concerns polygons with infinite number of vertices and provides various
distortion estimates in terms of intrinsic geometric characteristics of polygons.

In particular, this implies the coarse upper and lower estimates for the Grunsky and Teichm\"{u}ller norms of a
conformal map of the disk onto any piecewise $C^{1+}$-smooth  bounded quasicircle.
\end{abstract}

\date{\today\hskip4mm ({\tt PolygLine(3).tex})}

\maketitle

\bigskip

{\small {\textbf {2010 Mathematics Subject Classification:} Primary: 30C55, 30C62, 30F60; Secondary 30C80, 32G15, 46G20}

\medskip

}

\bigskip

\markboth{S. L. Krushkal}{Quantitaive theory of reflections across quasiconformal polygonal lines}\pagestyle{headings}

\bigskip\bigskip
\centerline{\bf 1. INTRODUCTORY REMARKS AND KNOWN RESULTS}

\bigskip\noindent
{\bf 1.1. General remarks}.
One of the important problems of geometric and quasiconformal analysis, also arising in applications, is
the quantitative evaluation of basic quasiinvariants of Jordan curves and arcs, such as the minimal dilatations
of quasiconformal continuations and reflections across these curves, Fredholm eigenvalues, etc.
This problem still has not been solved completely even for polygons.

This paper continues the author's research in quantitative investigation of this problem. It can be also regarded as
a complement to surveys \cite{Kr5}, \cite{Kr8}.

\bigskip\noindent
{\bf 1.2. Known results for polygons}.
The most general results established in this direction are the following two theorems.

\bigskip \noindent
{\bf Theorem 1}. \cite{Kr4}, \cite{Kr9} {\it For every unbounded convex domain $D \subset \C$ with piecewise $C^{1+\dl}$-smooth
boundary $L \ni \iy \ (\dl > 0)$ formed by a finite or countably many number of quasiintervals  and for all its fractional linear
images, we have the equalities
 \be\label{1}
\vk(f) = \vk(f^*) = k(f) = k(f^*) =  q_L = 1/\rho_L = 1 - |\a|,
\end{equation}
where $f$ and $f^*$ denote the appropriately normalized conformal maps $\D \to D$ and $\D^* \to D^* = \hC \setminus \overline{D}$,
respectively; $\pi \a$ is the opening of the least interior angle between the boundary arcs $L_j \subset L$;
$q_L$ and $\rho_L$ denote, respectively, the reflection coefficient and the Fredholm eigenvalue of
the boundary $L = \partial D$. Here $0 < \a < 1$ if the corresponding vertex is finite and $- 1 < \a < 0$ for the angle at the
vertex at infinity.

The same is true also for the unbounded concave domains (the complements of convex ones) which do not contain the infinite point; for those one must replace the last term by $|\beta| - 1$, where $\pi |\beta|$ is the opening of the largest interior angle of $D$. }

\bigskip
Here $\D$ denotes the unit disk $\{|z| < 1\}$ and  $\D^*$ its complement $\{z \in \hC = \C \cup \{\iy\}: \ |z| > 1\}$.
In the case of infinite number of vertices and sides of $D$, the endpoints of its boundary quasiintervals must accumulate only
to a countable set on $L$, and the modulus $\pi |\a_0|$ of opening of the angle at the limit vertex $A_0$ must be understand as the maximum of openings of angles adjoined to $A_0$ from inside of $D$  and for which the intervals $(A_j, A_0)$ are one of their sides.

The assumption of convexity in Theorem 1 is crucial and cannot be omitted for generic curvelinear polygons. A remarkable fact is
that in the case of the rectilinear polygons with a finite number of vertices the convexity is not needed.

Let a rectilinear polygon $P_n$ have the finite vertices $A_1, A_2, ... \ , A_{n-1}$ and the vertex $A_\iy = \iy$,
and let the interior angle at the vertex $A_j$ be equal to $\pi \a_j$  and be equal to $\pi \a_\iy$
at $A_\iy$, where $\a_\iy < 0$ and all $a_j \ne 1$, so that
$$
\a_1 + \dots + \a_{n-1} + \a_\iy = 2.
$$
Denoting the conformal map of the upper half-plane $U = \{z: \ \Im z > 0\}$ onto $P_n$ by $f_n$, we have

\bigskip\noindent
{\bf Theorem 2}. \cite{Kr7} {\it For any rectilinear polygon $P_n$, its conformal mapping function
$f = f_n \circ \sigma: \ \D \to P_n$ satisfies
 \be\label{2}
\vk(f) = k(f) = q_{\partial P_n} = 1/\rho_{\partial P_n} = |1 - |\a||.
\end{equation}
Here $\sigma$ is appropriate Moebius map $\D \to U$ and}
 \be\label{3}
|1 - |\a|| = \max \ \{|1 - |\a_1||, \dots , |1 - |\a_{n-1}||, \ |1 - |\a_\iy||\}.
\end{equation}

\bigskip
Such a result fails for the circular polygons, even for circular quadrilaterals.

Let us mention also that the proof of Theorem 1 relays on the properties of Finsler metrics on the negative generalized Gaussian curvature, while the proof of Theorem 2 involves the features of holomorphic motion generated by the pre-Schwarzian derivative
$b_f = f^{\prime\prime}/f^\prime$ of the Christoffel-Schwarz integral and of its extension.

\bigskip\bigskip
\centerline{\bf 2. INFINITE RECTILINEAR POLYGONS}

\bigskip
Our first aim is to find an extension of Theorem 2 to polygons with infinite number of vertices. In contrast to Theorem 1, 
there appear obstructions caused by the rigidity of rectilinearity.

For simplicity, we consider the case when the vertices of polygon accumulate to a countable set of limit points. The openings
of angles with vertices at the accumulation points are defined similar to Theorem 1.

Assume additionally that there is a subsequence $\{n_j\} \subset \mathbb N$ of vertices $A_n$, such that each $A_{n_j}$
can be joint with the infinite point by a ray $[A_{n_j}, \iy]$ inside the polygon $P$ forming $(n_1 + 1)$-gon $\wt P_{n_j+1}$,
and let the angle at $A_{n_j}$ with sides $[0, A_{n_j}]$ and $[A_{n_j}, \iy]$ in $\wt P_{n_j+1}$ does not exceed the deviation
 \be\label{4}
\beta = \sup_n |1- \a_n| < 1, \quad |1 - |\a_\iy||
\end{equation}

The first main result of this paper is

\bigskip\noindent
{\bf Theorem 3}. {\it Let unbounded rectilinear polygon $P$ have a countable set of finite vertices $A_1, A_2, \dots , A_n, \dots$ accumulating to $A_\iy = \iy$, and let its angles $\pi \a_n$ satisfy
$$
0 < \sup_n |1- \a_n| < 1, \quad |1 - |\a_\iy|| < 1.
$$
Then the conformal map $f(z)$ of the unit disk onto $P$ has the equal Grunsky and Teichm\"{u}ller norms, given by
  \be\label{5}
\vk(f) = k(f) = q_{\partial P} = 1/\rho_{\partial P} = \max (\sup_n |1 - \a_n|, |1 - |a_\iy||).
\end{equation}
}

A more general case of polygons whose vertices accumulate to a nonwhere dence countable set of limit points is investigated in
a similar way.

The applications of the above theorems are presented in the last sections.

\bigskip\bigskip
\centerline{\bf 3. PROOF OF THEOREM 3}

\bigskip
We start with a conformal map $g(\z)$ of the upper half-plane $U = \{z: \ \Im z > 0\}$ onto $P$ and take a point
$i b$ with $b > 0$ such that at this point the Schwarzian
$$
S_g(i b) \ne 0.
$$
Then the function
$$
f(z) = g \circ \sigma^{-1}_b(z) = z + a_2 z^2 + a_3 z^3 + \dots,
$$
where
$$
\sigma_b(\zeta) = (\zeta - i b)/(\zeta + i b),
$$
maps the unit disk onto the given polygon $P$ (with $f(0) = g(ib)$), and
 \be\label{6}
S_f(0) = 6(a_3 - a_2^2) \ne 0.
\end{equation}
Passing to inversions $F(z) = 1/f(1/z) = z + b_0 + b_1 z^{-1} + \dots$ univalent on the complementary disk
$\D^*$, one can rewrite (6) in the equivalent form $S_F(\iy) = b_1 \ne 0$.

The assumption (6) does not affect the geometric characteristics of $P$; thus we can assume that the initial function $f$
in the statement of Theorem 3 satisfies this inequality.

Now take the indicated sequence of vertices $A_{n_j} = f(z_{n_j}), \ z_{n_j} \in \mathbb S^1$, convergent to the infinite point
$A_\iy = \iy$.
Every rectlinear $(n_j + 1)$-gon $\wt P_{n_j+1}$ with vertices $A_1, \dots, A_{n_j}, A_\iy$  $P_n$ satisfis the assumption of Theorem 1 and by this theorem its mapping functions $\wt f_{n_j}$ obey the relations (1).
As $j \to \iy$, these polygons exhaust initial polygon $P_L$ from inside.

Now we investigate what occurs in the limit as $j \to \iy$. Some arguments applied here are similar to the case of convex quasiconformal polygons investigated in \cite{Kr9}; so some common facts will be shortened.

The general properties of univalent functions and of $k$-quasiconformal maps with $k < 1$ provide only the locally uniform convergence of (quasiconformall extended) functions $\wt f_{n_j}$ to $f$  on $\C$, and the semicontinuity of $k(f)$ and $\vk(f)$ on $\Belt(D^*)_1$ and $\T$, yields the inequalities
$$
k(f) \le \lim\limits_{j \to \iy} k(\wt f_{n_j}), \quad \vk(f) \le \liminf\limits_{j \to \iy} \vk(\wt f_{n_j}).
$$
We have to establish that the assumptions of Theorem 3 provide equalities in both of these relations.

To prove this, we use the Grunsky coefficients $\a_{m n}(f^\mu)$ of function $f \in S$, which generate
the holomorphic maps
 \be\label{7}
h_{\x}(f^\mu) = \sum\limits_{m,n=1}^\iy \ \a_{m n} (f^\mu) x_m x_n : \ \Belt(D^*)_1 \to \D
\end{equation}
with fixed $\x = (x_n) \in S(l^2)$. Here $f^\mu$ denotes the quasiconformal extension of $f$ to $D^*$, and $\mu$ is its Beltrami coefficient. Note that
 \be\label{8}
\sup_{\x \in S(l^2)} |h_{\x}(f^\mu)| = \vk(f^\mu).
\end{equation}

Another important tool is given by the homotopy functions
 \be\label{9}
f_t(z) = f(z, t) = \fc{1}{t} f(t z) = z + a_2 t z^2 + a_3 t^2 z^3 + \dots : \ \D \times \D \to \C.
\end{equation}
This complex homotopy is a special case of holomorphic motions. It admits the following needed for us properties, presented
by the next three lemmas.

\bigskip\noindent
{\bf Lemma 1}. {\it For any homotopy function $f_t$, the following equalities are valid: }
$$
\lim\limits_{|t| \to 1} k(f_t) = k(f), \quad \lim\limits_{|t| \to 1} \vk(f_t) = \vk(f).
$$

\bigskip\noindent
{\bf Lemma 2} {\it For any function $f \in S$ satisfying (6), the extremal
quasiconformal extensions of the homotopy functions $f_t$ to the disk $\D^*$ are defined for sufficiently small $|t| \le r_0 = r_0(f) \ (r_0 > 0)$ by nonvanishing holomorphic quadratic differentials, and therefore,
$\vk(f_t) = k(f_t)$}.

\bigskip\noindent
{\bf Lemma 3}. {\it If the homotopy function $f_t$ of $f \in S$ satisfies $\vk(f_{t_0}) = k(f_{t_0})$ for
some $0 < t_0 < 1$, then the equality $\vk(f_t) = k(f_t)$ holds for all $|t| \le t_0$ and the homotopy disk
$\D_{f} = \{S_{f_t}\} \subset \T$ has no critical points $t$ with $0 < |t| < t_0$.
}

\bigskip
The proofs of these lemmas are given in \cite{Kr8}, \cite{Kr9}, \cite{Ku2}, \cite{Ku3}).
Note that if a sequence $\{f_n\} \subset S$ is convergent to $f$ locally uniformly on $\D$, then, for every
$|t| < 1$, the Schwarzians $S_{f_{n,t}}$ of the homotopy functions are convergent to $S_{f_t}$ in the norm of the space $\B$.

\bigskip
In view of the equality (8), one can select the sequences of small $\ve_n \searrow 0, \ t_n \nearrow 1$ and of points $\x^{(n)} \in S(l^2)$ such that
 \be\label{10}
|h_{\x^{(n)}}(f_{t_n}^{\mu_{t_n}})| \ge \vk(f) - \ve_n, \quad n = 1, 2, \dots;
\end{equation}
here $\mu_t$ denote the extremal Beltrami coefficients of $f_t$ in $\D$ (all these $\mu_t$ are of Teichm\"{u}ller type).

Now, using the functions $h_{\x^{(n)}}(f_t^{\mu_t}): \ \D \to \D$, we pull back the hyperbolic metric $\ld_\D(t)|dt| = |dt|/(1 - |t|^2)$ of the disk $\D$ onto this disk , getting conformal metrics $\ld_{h_\x^{(n)}}(t) |dt|$ on $\D$ with
$$
\ld_{h_{\x^{(n)}}}(t) = \fc{|h_{\x^{(n)}}^\prime (t)|}{(1 - |h_{\x^{(n)}}(t)|^2)}
$$
of Gaussian curvature $-4$ at noncrical points.
We take the upper envelope
$$
\ld_\vk(t) =  \sup_n \ \ld_{h_{\x^{(n)}}}(t)
$$
of these metrics followed by its upper semicontinuous regularization. This provides a logarithmically subharmonic metric
on the unit disk, for which we preserve the notation $\ld_\vk(t)$.

On a standard way, one obtains that $\ld_\vk$ has at any its noncritical point $t_0$  a supporting subharmonic metric $\ld_0$  (i.e., such that $\ld_0(t_0) = \ld(t_0)$ and $\ld_0(t) \le \ld(t)$ in a neighborhood of $t_0$) of Gaussian curvature at
most $- 4$, and hence, $\kappa_{\ld_\vk} \le - 4$ (the details see, e.g., in \cite{Ah1}, \cite{He}, \cite{Kr3}).

On the other hand, Lemma 3 and the inequality (10) imply that the constructed metric $\ld_\vk(t)$ must be equal at the origin $t = 0$  to the infinitesimal Kobayasi-Teichm\"{u}ller metric $\ld_{\mathcal K}(t)$:
 \be\label{11}
\ld_{\vk}(0) = \ld_{\mathcal K}(0).
\end{equation}

On the homotopy disk $\{|t| < 1\}$, we have the equality $\ld_{\vk}(t) \le \ld_{\mathcal K}(t)$ following from the relation between the Grunsky and Teichm\"{u}ller norms of $f$, and our goal now is to establish that in fact these metrics must be
equal on this disk.
We shall use their  (generalized) Gaussian curvatures defined for the upper semicontinuous Finsler metrics
$ds = \ld|dt|$ in a domain $\Om \subset \C$ by
$$
\kappa_\ld (t) = - \fc{\Delta \log \ld(t)}{\ld(t)^2},
$$
where $\Delta$ is the generalized Laplacian
$$
\Delta \ld(t) = 4 \liminf\limits_{r \to 0} \frac{1}{r^2} \Big\{ \frac{1}{2 \pi} \int_0^{2\pi} \ld(t + re^{i \theta})
d \theta - \ld(t) \Big\}
$$
(provided that $- \iy \le \ld(t) < \iy$). In fact, this is equivalent to regard the differential operator
$\Delta = 4 \partial^2/\partial z \partial \ov{z}$ in the distributional sense.

Similar to $C^2$ functions, for which $\Delta$ coincides with the usual Laplacian, one obtains that $\ld$ is
subharmonic on $\Om$  if and only if $\Delta \ld(t) \ge 0$; hence, at the points $t_0$ of local maximuma of $\ld$
with $\ld(t_0) > - \iy$, we have $\Delta \ld(t_0) \le 0$.

An important fact is that the metrics $\ld$ satisfying the inequality
$$
\Delta \log \ld \ge K \ld
$$
with $ K = \const > 0$
(then $u = \log \ld$ can be negative) admit the following {\bf Minda's maximum principle} \cite{Mi}, which is a deep
variation of the Ahlfors-Schwarz lemma:

\bigskip\noindent
{\bf Lemma 4}. {\it If a function $u : \ D \to [- \iy, + \iy)$ is upper semicontinuous in a domain
$D \subset \C$ and its generalized Laplacian satisfies the inequality $\Delta u(z) \ge K u(z)$ with some positive
constant $K$ at any point $z \in D$, where $u(z) > - \iy$, and if
$$
\limsup\limits_{z \to \z} u(z) \le 0 \ \ \text{for all} \ \z \in \partial D,
$$
then either $u(z) < 0$ for all $z \in D$ or $u(z) \equiv 0$ on $D$.  }

\bigskip
Now choose a sufficiently small neighborhood $U_0$ of the origin $t = 0$ and put
$$
M = \{\sup \ld_{\mathcal K}: t \in U_0\}.
$$
Then in this neighborhood, we have $\ld_{\mathcal K}(t) + \ld_\vk(t) \le 2M$. Consider the ratio
$$
u = \log \fc{\ld_\vk}{\ld_{\mathcal K}}.
$$

Then (cf. \cite{Mi}, \cite{Di}) for $t \in U_0$,
$$
\Delta u(t) = \Delta \log \ld_\vk(t) - \Delta \log \ld_{\mathcal K}(t)
\ge 4 (\ld_\vk^2 - \ld_{\mathcal K}^2) \ge 8M (\ld_\vk - \ld_{\mathcal K}),
$$
and the elementary estimate
$$
M \log(t/s) \ge t - s \quad \text{for} \ \ 0 < s \le t < M
$$
(with equality only for $t = s$) implies that
$$
M \log \fc{\ld_\a(t)}{\ld_\vk(t)} \ge \ld_\a(t) - \ld_\vk(t);
$$
hence, $\Delta u(t) \ge 8 M^2 u(t)$.

Noting that by (10) $u(0) = 0$ and applying to $u$ Lemma 4 with $K = M^2$ on the neighborhood $\Om = U_0$,
one derives that both metrics $\ld_\vk$ and $\ld_{\mathcal K}$ must be equal on $U_0$,
and in the similar way their equality on the entire disk $\D$. This proves the infinitesimal version of Theorem 3.

To get its global version (4), we apply the following reconstruction lemma for the Grunsky norm proven in \cite{Kr3}, which provides that this norm is the integrated form of $\ld_\vk$ along the Teichm\"{u}ller extremal disks.

\bigskip\noindent
{\bf Lemma 5}. {\it On any  Teichm\"{u}ller extremal disk $\D(\mu_0) = \{t \mu_0/\|\mu_0\|_\iy: |t| < 1\}  \subset \Belt(\D^*)_1$,
we have the equality}
$$
\tanh^{-1}[\vk(f^{r\mu_0/\|\mu_0\|_\iy})] = \int\limits_0^r \ld_\vk(t) dt.
$$

Integrating the metrics  $\ld_\vk$ and $\ld_{\mathcal K}$ along the indicated extremal disks, one derives from the above
relations the desired equality (4), completing the proof of Theorem 3.

\bigskip\bigskip
\centerline{\bf 4. EXTENSION AND APPLICATIONS OF THEOREMS 1-3.}
\centerline{\bf ADDITIONAL REMARKS}

\bigskip\noindent
{\bf 4.1. Extensions of Theorem 3}.
Though rectilinearity causes rather strong rigidity for extensions of Theorem 3, the arguments applied in its proof can be extended to some more general domains whose conformal mapping functions appear as the limit functions  of sequences of conformal maps of $\D$ onto polygons under locally uniform convergence on the disk $\D$. This provides, in particular, the domains whose
boundaries contain uncountable limit sets for vertices $A_n$, and Theorem 3 yields explicitly the dilatations for the limit maps.
In view of the indicated rigidity, such extensions are very restricted.

\bigskip\noindent
{\bf 4.2. Starlike polygons}. The hypotheses of Theorem 3 are fulfilled for the infinitely edged rectilinear
polygons which are starlike with respect to some their point $z_0$, that means, for any point $z^* \in \ov P$ the segment
$[z_0, z^*]$ is placed entirely in $\ov P$ and its part in $P$ is a connected half-interval.

Without loss of generality, one can assume that $z_0 = 0$, hence deal with the intervals of rays $tz, \ t > 0$.
For all such polygons we have

\bigskip\noindent
{\bf Theorem 4}. {\it Let $P$ be a starlike rectilinear polygon $P$ with at most countable set of finite vertices $A_1, A_2, \dots , A_n, \dots$ and a vertex $A_\iy = \iy$. If its deviation (4) satisfies $0 < \beta < 1$ and the conformal mapping function
$f: \ \D \to P$ obeys the assumption (6), then the associated quasiinvariats of $P$ are given by
 \be\label{12}
\vk(f) = k(f) = q_{\partial P} = 1/\rho_{\partial P} = \beta.
\end{equation}
}
\bigskip\noindent
{\bf 4.3. Quasiconformal distortion of bounded polygons}.
As was indicated above, the assumption that one of vertices is placed at the infinite point is crucial. The assertions similar
to Theorems 1-3 fail in the case of bounded polygons (this also is illustrated by example {\bf 5.3} below).

However these theorems allow one to establish for bounded polygons a simple upper bound estimating the possible distortions. This
coarse bound also estimates the associated quasiinvariants in terms of intrinsic geometric characteristics of a given polygon.

To formulate the results, we shall need some preliminary construction.
Let $\mathcal A$ be a given angle of opening $\a$ with vertex $A \in \C$, whose sides are  two rays $L_1, \ L_2$ outgoing from $A$
to $\iy$. Take its bisector $L_0$ outgoing from $A$ and dividing $\mathcal A$ on two equal parts.

For each point $t \in L_0$, we construct the rhombus $\mathcal R_t$ with equal angles at the opposite points $t$ and $A$.
Two from its sides are placed on the rays $L_t^{(1)}, \ L_t^{(2)}$ outgoing from this point to infinity. These rays intersect the lateral sides $L_1, \ L_2$
of the original angle $\mathcal A$ in the points $B_{1t}, \B_{2 t}$, which are the additional vertices of rhombus $\mathcal R_t$.
The sides of different $\mathcal R_t$ are pair-wise parallel, their direction and the angles of rhombuses do not depend on $t$.
We regard each $\mathcal R_t$ as a {\it rhombus associated with} $\mathcal A$.

Consider also the angle $\a_{1t}$ with vertex at the point $B_{1t}$ and the sides along $L_1, \ L_t^{(1)}$ and the similar angle
$\a_{2t}$ with vertex at $B_{2t}$. In view of parallelism, also these angles are independent on $t$. Denote their common values
by $\a_{1, A}, \ \a_{2, A}$ and call these quantities the {\it  adjoint rhombic angles for vertex} $A$.

Now we may formulate the distortion theorems for bounded polygons.
We start with rectilinear polygons $P_n$ with a finite number of vertices and assume that all angles $\pi \a_j < \pi$, i.e., their
factors $\a_j < 1$.

\bigskip\noindent
{\bf Theorem 5}. {\it The Teichm\"{u}ller and Grunsky norms of the inner conformal mapping function $f$ of any bounded rectilinear polygon $P$ with finite or a countable set of vertices $A_n \in \C$ with angles $\pi \a_n$ are estimated by
  \be\label{13}
\max_j |1 - \a_j| \le \vk(f) \le k(f) = q_{\partial P} < \max_j |1 - \a_j| + b_P < 1.
\end{equation}
where
$$
b_P = \max(1 - \a_{1,0}^{ad}, 1 - \a_{2,0}^{ad})
$$
and $\pi a_{1,0}^{ad}, \ \pi \a_{2,0}^{ad}$ are the adjoint rhombic angles for vertex $A_{j_0}$ of one of a smallest inner angle of $P$.  }

\bigskip\noindent
{\bf Prof}. The left-hand side inequality in (13) is a consequence of the following lemma.

\bigskip\noindent
{\bf Lemma 6}. {\it If a closed curve $L \subset \hC$ contains two $C^{1+}$-smooth arcs with the interior intersection
angle $\pi \a$, then its least nontrivial Fredholm eigenvalue}
$$
\fc{1}{\rho_L} \ge |1 - |\a||.
$$

This remarkable fact was established for angles with analytic sides by K\"{u}hnau \cite{Ku1}.
It extends to arbitrary $C^{1+\delta}$-smooth sides by appropriate approximation.

As for other inequalities in (13), consider first the polygons $P$ with finite number of vertices. One can assume that the
minimal angle of $P$ is $A_1$.
Applying to this angle the above construction, one immediately obtains the inequalities (13) as a consequence of Theorem 2.

Then the case of infinite polygons is investigated in similar way applying some arguments from the proof of Theorem 3.

The estimate (13) trivially yields the upper and lower bounds for dilatation of extremal reflections across $\partial P$ and
for Fredholm eigenvalue $\rho_{\partial P}$; these bounds also depend only on angles of $P$.

\bigskip\noindent
{\bf 4.4}. If the mapping function $f$ satisfies (6), then in (13) $\vk(f) = k(f)$.

Indeed, under this assumption the univalent functions $f_n$ with $\vk(f_n) = k(f_n)$ cannot approximate the maps $f$ with
$\vk(f) < k(f)$ (even locally uniformly). This is obtained by applying the homotopy functions $f_t(z)$ (cf. \cite{Kr6}).

\bigskip\noindent
{\bf 4.5}.
In the case when the polygon $P$ has a symmetry axes passing through the distinguished vertex $A_{J_0}$, the adjoint angles $\pi \a_{1,j_0}^{ad}$ and $\a_{2,j_0}^{ad}$ are equal. Then
$$
b_P = \max(1 - \a_{1,0}^{ad})
$$
and this angle $\pi \a_{1,0}^{ad}$ is simply determined by $\pi \a_1$.  This case has an intrinsic interest.

\bigskip\noindent
{\bf 4.6}. Applying approximation, one can derive from Theorem 5 a similar result for arbitrary piece-wise smooth quasiconformal
domains.

\bigskip\noindent
{\bf Theorem 6}. {\it The conformal map of the disk onto any piece-wise smooth bounded quasicircle $L \subset \C$ with finite number of singular points obeys the deformation bounds (13) determined by the smallest angle between the smooth intervals
of $L$. }

\bigskip
The existence of singular points is here necessary. The assumptions on regularity of the boundary $L$ in
Theorem 6 can be essentially weakened.

In particular, all this holds for circular polygons.

\bigskip\bigskip
\centerline{\bf 5. EXAMPLES}

\bigskip\noindent
{\bf 5.1. Example 1: rectangle}.
It follows from Theorem 1 that for any (oriented) closed unbounded curve $L$  with the convex interior which is $C^{1+\dl}$ smooth at all finite points and has at infinity the asymptotes approaching the interior angle $\pi \a_0 <0$, we have the exact representation
of the associated distortions:
$$
\vk(f) = k(f) = q_L = 1/\rho_L = 1 - |\a_\iy|,
$$
where $f$ maps conformally the disk $\D$ onto the inner domain of $L$.

As for bounded polygons, we do not have such exact values, but Theorems 5 and 6 yield the following coarse bounds.

Let $D$ be a bounded convex domain with one angular boundary point $z_0$, and let the $L = \partial D$ be smooth in all its
other points. Then Theorem 6 implies
$$
k(f) < \a_0 + b_D,
$$
where $\pi \a_0$ is the opening of the angle with vertex $z_0$ and $b_0$ is determined by the limit of associated rombic angles.

Of course, the additional term can be improved or even omitted. For example, in the case of rectangles $P_4$ with modulus $\mu$
sufficiently close to $1$, we have the value
 \be\label{14}
k(f) = 1/2
\end{equation}

and $k(f) = q_{\partial P_4} > 1/2$ for the rectangles with modulus $\mu > 2.76$ (see \cite{We}).

In contrast, by \cite{Kr4} for any rectangle the equality $\vk(f) = k(f)$ is valid.

Rounding the rectangle $P_4$ off near three of its vertices, one obtains a convex domain $D$ with one angular boundary point,
for which the estimate (14) is preserved.

\bigskip\noindent
{\bf 5.2. Example 2: infinite ladders}.
Let a quasicircle $L$ contains both points $0$ and $\iy$, and is a union of countable number of mutually orthogonal crossbars
so that near $z = 0$ it is an interval of the real axes. Assume that the steps of $L$ increase unlimitedly by approach to $\iy$
(such lines appear often in numerical analysis). This quasicicle divides the plane onto two polygons with infinite number
of vertices, and their conformal mapping functions are quantitatively characterized by the equalities (5).
The case of a quasiinterval $L$ is investigated by Theorem 5.
The same is valid for images $\gamma(L)$ under the Moebius maps.

The invariants of quasicircles with triangular, trapezoidal and of other forms teeth are estimated in a similar way

\bigskip\noindent
{\bf 5.3. Example 3: the Koch snowflake curves}.
Theorems 1-3 can be applied to coarse estimating the dilatations of rather pathological quasicircles. We illustrate this
on snowflake curves (fractal curves possessing selfsimilarity).
Such curves are broadly applied in complex dynamics and other fields of mathematics.

We present the snow-flake curves applying the construction given by Astala \cite{As}.
First recall that a quasiinterval $I$ is calles {\bf selfsimilar} if there are $N$ similarities
$\sigma_j, \ 1 \le j \le N$ and $N \ge 2$, such that
$$
I = \bigcup\limits_1^N \sigma_j(I)
$$
and for $j_1 \ne j_2$ the intersection $\sigma_{j_1}(I) \cap \sigma_{j_2}(I)$ is either empty or consists of a point.

Fix the points $z_1 = 0, \ z_2 = 2, \ z_3 = 3 + \sqrt{3}, \ z_4 = 4, \ z_5 = 6$ and denote by $\sigma_j$ the similarity which
maps the line segment $[z_1, z_5]$ to $[z_j, z_{j+1}]$ with $\sigma_j(z_1) = z_j$.
Letting $\sigma(E) = \bigcup_1^4 \sigma_j(E), \ E \subset \C$, one obtains that the iterated arcs
$\sigma^p[z_1, z_5] = \sigma(\sigma^{p-1}[z_1, z_5])$ approach as $p \to \iy$  in the Hausdorff metric
$$
d_H(E_1, E_2) = \max \Bigl\{ \sup_{x \in E_1} d(x, E_2), \sup_{y \in E_2} d(y, E_1)\Bigr\}
$$
a snowflake curve $\mathcal S_{1/3}$; here $d$ is the Euclidean metric (see \cite{Ha}).

The ratio $1/3$ of decreasing the length of the segment $[z_1, z_5]$ can be replaced by any $t$ from $(1/4, 1/2)$,
which provides the snowflake curve $\mathcal S_t$ of Hausdorff dimension $(\log \fc{1}{4})/\log t$ \cite{Ha}.

Denote by $l$ the polygonal interval joining the initial points $z_1, z_2, \dots, z_5$.
Noting that all iterations $\sigma^p$ preserve its endpoints $z_1, z_5$, we extend all images $\sigma^p(l)$ by translations
$z \mapsto z + 6 n, \ n = 0, \pm 1, \pm 2, \dots$, getting the unbounded polygonal lines and in the limit as $p \to \iy$
an unbounded snowflake quasicircle  $\mathcal S_t^\iy$.

The above theorems imply the exact values of any iterated curve and their upper and lower bounds of type (13) (without the term $b_P$) for the limit curve $\mathcal S_t^\iy$.

\bigskip\bigskip
\centerline{\bf 6. ESTIMATING QUASIREFLECTIONS ACROSS ARCS AND}
\centerline{\bf ARBITRARY SETS}

\bigskip\noindent
{\bf 6.1. General approach}.
One of the open and very intriguing problems in quasiconformal analysis and its applications is to estimate the reflection coefficients $q_E$ of arbitrary quasiconformal mirrors $E \subset \hC$ (sets remaining fixed under orientation reversing  quasiconformal involutions of the sphere). In contrast to the case of (closed) curves, only a few estimates have been established for quasireflections across the arcs and arbitrary sets.

The above theorems provide simultaneously the estimates for sets placed on (sufficiently regular) quasicircles with
distinguished points.
This includes, in particular, all finite collections of points on the extended complex plane $\hC$ which can be joined by
a polygonal line without self-intersection. An underlying fact here is  given by the following

\bigskip\noindent
{\bf Lemma 7}. {\it For any set $E \subset \hC$ which admits quasireflections, there is a quasicircle $L \supset E$ with the
same reflection coefficient; therefore,}
\be\label{15}
Q_E = \min \{Q_L : \ L \supset E \ \ \text{quasicircle}\}.
\end{equation}

\bigskip
This was established for finite sets $E = \{z_1, ... \ , z_n\}$ by K\"{u}hnau \cite{Ku2}, using Teichm\"{u}ller's theorem on extremal quasiconformal maps, and in the general case in \cite{Kr2}.

\bigskip
Lemma 7 and Theorem 1 provide a possibility to estimate the distortion of subsets of rather general quasicircles.

As an example, one can use any subset  of a polygonal line $L = \partial P$ containing the subarcs
forming the smallest angle $\pi \a$. This line is an extremal quasicircle for $E$, on which the minimum in (15) is attained.
Hence, for any subset $e \subset L \setminus E$ the union $E \cup e$ has the same reflection coeffficient $1 - |\a|$
as the initial set $E$.

\bigskip\noindent
{\bf 6.2. Reflections across analytic arcs}.
In the more restrictive case of analytic arcs, one can also apply another approach to estimating the reflection dilatations.
This approach was developed in \cite{Kr1}, \cite{Kr5} and relates to the best approximation of holomorphic functions
estimated by the pluricomplex Green function of compact sets in $\C^n$. The estimates of the reflection dilatations obtained
on this way are sharp.

First recall that this function with pole at infinity is defined for a given compact $k \subset \C^n$ by
defined by
$$
g_K(z) = \sup \{u(z) : \ u \vert K \leq 0 \}
$$
(and followed upper semicontinuous regularization of this function), where the supremum is taken over all plurisubharmonic functions $u$ on $\C^n$ with growth $u(z) = \log |z| + O(1)$ as $z \to \iy$  .
In particular,
$$
g_{[-1, 1]^n} (z) = \max\limits_{1 \le j \le n} \ \{\log |h (z_j)|: \ z = (z_1, ..., z_n) \in \C^n \} ,
$$
where $h(t) = t + \sqrt{t^2 - 1}$ is the single-valued branch of inverse function to the Joukovski function chosen so
that $h(t) > 1$ for $t > 1$.

Due to the classical Bernstein-Walsh-Ciciak theorem, a continuous function $h$ on a compact $K \subset C^n$ with
$g_K(z)|K \equiv 0$ extends holomorphically to the region
$$
G_R = \{ z \in {\mathbb C}^n : \ g_K (z) < \log R \}, \quad R > 1,
$$
and $R$ satisfying
$\limsup\limits_{m \to \iy} e_m^{1/m} (h, K) = 1/R, \quad R > 1$.

We apply it to the function
$$
F_f(x, \xi) = \log \frac{f(x) - f(\xi)}{x - \xi} : \ [-1, 1]^2 \to \C,
$$
where the single-valued branch of $\log w$ is chosen to be positive for $w > 1$, and $F(x, x) = \log f'(x)$.
This function is holomorphic on $[-1, 1]^2$ for injective $f$.

Now, given an analytic arc $L = g(l)$, where $l$ is a subarc of the unit circle and $g$ is injective on
$\gamma$, then applying Theorem 3.2 from \cite{Kr5}, one derives that {\it the reflection coefficient of $L$ for any injective analytic function $g$ on the boundary arc $l$ is estimated by
 \be\label{16}
{r^2 + 1/r^2} = \frac{1}{a^2 + b^2},
\end{equation}
with
$$
r = \limsup\limits_{m \to \iy} \ e_m^{1/m} \bigl(F_{g \circ \gamma}, l \times l \bigr) < 1,
$$
where $\gamma$ is the Moebius transformation of $\hC$ moving the arc $l$ onto the segment $[-1, 1]$
and $a, b$ are the semiaxes of the ellipse with foci at -1, 1 such that $a + b = \fc{1}{r}$. The bound (16) is sharp}.

\medskip
{\small \em{ \leftline{Department of Mathematics, Bar-Ilan University, 5290002 Ramat-Gan, Israel}
\leftline{and Department of Mathematics, University of Virginia, Charlottesville, VA 22904-4137, USA}
}}


\bigskip
\bigskip
\begin{thebibliography}{EKK}
{\small

\bibitem{Ah1}
L. Ahlfors, {\it An extension of Schwarz's lemma}, Trans. Amer. Math. Soc. \textbf{43} (1938), 359-364.

\bibitem{Ah2}
L.V. Ahlfors, {\it Lectures on Quasiconformal Mappings}, Van Nostrand, Princeton, 1966.

\bibitem{As}
K. Astala, {\it Selfsimilar zippers}, Holomorphic Functions and Moduli I, Math. Sci. Research Inst. Publications, vol. 10
(D. Drasin et al., eds.), New York, 1988, 61 - 73.

\bibitem{Di}
S. Dineen, {\it The Schwarz Lemma}, Clarendon Press, Oxford, 1989.

\bibitem {GL}
F.P. Gardiner and N. Lakic, {\it Quasiconformal Teichm\"{u}ller Theory}, Amer. Math. Soc., Providence, RI, 2000.

\bibitem{Gr}
H. Grunsky, {\it Koeffizientenbedingungen f\"{u}r schlicht abbildende meromorphe Funktionen},
Math. Z. \textbf{45} (1939), 29-61.

\bibitem{Ha}
J.E. Hutchinson, {\it Fractals and selfsimilarity}, Indiana Univ. Math.J. \textbf{30} (1981), 713-747.

\bibitem{He}
M. Heins, {\it A class of conformal metrics}, Nagoya Math. J. \textbf{21} (1962), 1-60.

\bibitem{Kr1}
S.L. Krushkal, {\it On the best approximation and univalence of holomorphic functions},
Complex Analysis in Contemporary Mathematics. In honor of 80th Birthday of Boris Vladimirovich Shabat (E.M. Chirka, ed.),
Fasis, 2001, pp. 153-166.

\bibitem{Kr2}
S. L. Krushkal,
{\it Quasiconformal reflections across arbitrary planar sets}, Scientia, Series A: Mathematical Sciences \textbf{8} (2002),
57-62.

\bibitem{Kr3}
S.L. Krushkal, {\it Plurisubharmonic features of the Teichm\"{u}ller  metric}, Publications de
l'Institut Math\'{e}matique-Beograd, Nouvelle s\'{e}rie \textbf{75(89)} (2004), 119-138.

\bibitem{Kr4}
S.L. Krushkal, {\it Quasireflections, Fredholm eigenvalues and Finsler metrics}, Doklady Mathematics \textbf{69}
(2004), 221-224.

\bibitem{Kr5}
S.L. Krushkal, {\it Quasiconformal extensions and reflections}, Ch. 11 in: Handbook of Complex Analysis: Geometric
Function Theory, Vol. II (R. K\"{u}hnau, ed.), Elsevier Science, Amsterdam, 2005, pp. 507-553.

\bibitem{Kr6}
S.L. Krushkal, {\em Strengthened Moser's conjecture, geometry of Grunsky coefficients and Fredholm eigenvalues},
Central European J. Math \textbf{5(3)} (2007), 551-580.

\bibitem{Kr7}
S.L. Krushkal, {\it On Fredholm eigenvalues of unbounded polygons}, Siberian Math. Journal  \textbf{60} (2019),
no. 5, 896-901.

\bibitem{Kr8}
S.L. Krushkal, {\it Fredholm eigenvalues and quasiconformal geometry of polygons}, J. Math. Sci. \textbf{252}(4) (2021),
472-501. DOI: 10.1007/s10958-020-05175-4

\bibitem{Kr9}
S.L. Krushkal, {\it Analytic and geometric quasiinvariants of convex curvelinear polygons with infinite number of
vertices}, J. Math. Sci. (2024), to appear.

\bibitem{Ku1}
R. K\"{u}hnau, {\it M\"{o}glichst konforme Spiegelung an einer Jordankurve}, Jber. Deutsch. Math. Verein.
\textbf{90} (1988), 90-109.

\bibitem{Ku2}
R. K\"{u}hnau, {\it Interpolation by extremal quasiconformal Jordan curves}, Siberian Math. J. \textbf{32} (1991),   257-264.

\bibitem{Ku3}
R. K\"{u}hnau, {\it Wann sind die Grunskyschen Koeffizientenbedingungen hinreichend f\"{u}r $Q$-quasikonfor\-me Fortsetzbarkeit}?
Comment. Math. Helv. \textbf{61} (1986), 290-307.

\bibitem{M}
I.M. Milin, {\it Univalent Functions and Orthonormal Systems}, Transl. of Mathematical Monographs, vol. 49,
Transl. of Odnolistnye funktcii i normirovannie systemy, Amer. Math. Soc., Providence, RI, 1977.

\bibitem{Mi}
D. Minda, {\it The strong form of Ahlfors' lemma}, Rocky Mountain J. Math., \textbf{17} (1987), 457-461.

\bibitem{Po}
Chr. Pommerenke, {\it Univalent Functions}, Vandenhoeck $\&$ Ruprecht, G\"{o}ttingen, 1975.

\bibitem{Sc}
M. Schiffer, {\it Fredholm eigenvalues and Grunsky matrices}, Ann. Polon. Math. \textbf{39} (1981), 149-164.

\bibitem{We}
S. Werner, {\it Spiegelungskoeffizient und Fredholmscher Eigenwert f\"{u}r gewisse Polygone},
Ann. Acad.Sci. Fenn. Ser. AI. Math., \textbf{22} (1997), 165-186.





}


\end{thebibliography}
\end{document}